\documentstyle[12pt]{article}
\parskip 0.05 truein\baselineskip 4pt\lineskip 4pt\def\a{\alpha}
\def\b{\beta}
\def \F{\hbox{$I\hskip -4pt F$}}\def\Z{\hbox{$Z\hskip -5.2pt Z$}}
\def\sZ{\hbox{$\sc Z\hskip -4.2pt Z$}}

\def\vs{\vspace*} \def\Sl{{\rm sl}}
\def\nob#1{\mbox{#1}}
\def\qed{\hfill \hfill \ifhmode\unskip\nobreak\fi\ifmmode\ifinner
\else\hskip5pt\fi\fi
\hbox{\hskip5pt\vrule width4pt height6pt depth1.5pt\hskip 1 pt}}
\def\d{\delta}\def\g{\gamma}\def\l{\lambda}
\def\si{\sigma}
\def\sc{\scriptstyle}\def\ssc{\scriptscriptstyle}
\def\cl{\centerline}
\def\ol{\overline}
\def\wt{\widetilde}\def\rar{\rightarrow}

\def\bs{\backslash}
\def\hs{\hspace*}\def\rb{\raisebox}\def\ra{\rangle}\def\la{\langle}
\def\ni{\noindent}\def\hi{\hangindent}\def\ha{\hangafter}
 \def\f{ \,\,\,\forall \,\,\,}
\def\Vir{\mbox{\bf\rm Vir}} \def\bvv{{\ol V}}
\begin{document}
\cl{GENERALIZED VIRASORO AND SUPER-VIRASORO ALGEBRAS } \cl{ AND
MODULES OF THE INTERMEDIATE SERIES } \vs{4pt}\cl{(Appeared in
{J.~Alg.} {\bf252} (2002), 1-19.)}
\par\vs{-7pt}\
\par
\centerline{Yucai Su$^{\ast}$ and Kaiming Zhao$^{\dag}$}
\par{\small\it
* Department of Mathematics, Shanghai Jiaotong University, Shanghai 200030, China.

\dag Institute of Mathematics, Academia Sinica, Beijing 100010,
China.}
\par\vs{-4pt}\
\par
\cl{\bf               \S 1. Introduction
}
\par
Recently, a number of new classes of
infinite-dimensional simple Lie algebras over
a Field of characteristic 0 were discovered by several authors
(see the references at the end of this paper).
Among those algebras, are the generalized Witt algebras.
The higher rank Virasoro algebras was introduced by
Patera and Zassenhaus [PZ], which are 1-dimensional universal central extensions
of some generalized Witt algebras, and the higher rank super-Virasoro algebras
was introduced by Su [S].
\par
In this paper, we will introduce and study generalized Virasoro and super-Virasoro algebras which are slightly more general than higher rank
Virasoro and super-Virasoro algebras. In Section 2, we give the definition and all the automorphisms of generalized
Virasoro algebras, which are different  from that given in [PZ] for higher rank
Virasoro algebras. In Section 3, we exhibit many two dimensional subalgebras of generalized
Virasoro algebras,  which should occur but
were not  collected in [PZ] for higher rank Virasoro algebras.
We believe our results are correct although they are different from those given in [PZ].  Section 4 devotes to
 determining all modules of the intermediate series over the
generalized Virasoro algebras, see Theorem 4.6.
Then in Section 5  we introduce the notion of
generalized super-Virasoro algebras and determine
the modules of the intermediate series over the
generalized super-Virasoro algebras, see Theorem 5.4. Our work here
generalizes the results
of [PZ] and [S] in the following two ways, first, the ground field is
any field of characteristic 0, not necessarily algebraically closed;
second, the subgroup $M$ of $\F$ is arbitrary, not necessarily of finite rank.
\par\
\par\
\par
\cl{\bf\S 2. Generalized Virasoro algebras}
\par
Let $A$ be a torsion free abelian additive group, $\F$ any field of characteristic 0, $T$ a vector space over $\F$.
Denote by $\F A$ the group algebra of $A$ over $\F$. The elements
$t^x,\, x\in A$, form a basis of this algebra, and the multiplication
is defined by $t^x\cdot t^y=t^{x+y}$. The tensor product
$W=\F A\otimes_{\F} T$ is a free left $\F A$-module. We denote an
arbitrary element of $T$ by $\partial$. We shall write $t^x\partial$
instead
of $t^x\otimes\partial$. Choose a pairing $\phi:T\times A\rar \F$ which
is
$\F$-linear in the first variable and additive in the second one. For
convenience we shall also use the following notations:
$$
\phi(\partial,x)=\la\partial,x\ra=\partial(x),\ \ \f \partial\in T,\ x\in
A.
\eqno(2.1)$$
The following bracket
$$
[t^x\partial_1,t^y\partial_2]=
t^{x+y}(\partial_1(y)\partial_2-\partial_2(x)\partial_1),
\ \ \f x,y\in A,\ \partial_1,\partial_2\in T,
\eqno(2.2)$$
makes $W=W(A,T,\phi)$ into a Lie algebra, which was referred to as a
{\it generalized Witt algebra} and studied in [DZ2]. The following result
was proved there.
\par\ni
{\bf Theorem 2.1}. {\it (1) $W=W(A,T,\phi)$ is simple if and only if
$A\ne0$ and $\phi$ is nondegenerate.
\vs{-3pt}\par
(2) If $\theta: W\rar W'$ is a  Lie algebra isomorphism between
two simple algebras $W=W(A,T,\phi)$ and $W'=W(A',T',\phi')$,
then
there exists $\chi\in{\rm Hom}(A,\F^*)$ where  $\F^*$ is the
\nob{multiplicative} group $\F\setminus\{0\}$,
isomorphisms $\si:A\rar A'$ and $\tau:T\rar T'$ satisfying
$\la\tau(\partial),\sigma(x)\ra=\la\partial,x\ra$ such that
$\theta(t^x\partial)=\chi(x)t^{\si(x)}\tau(\partial)$.
\vs{-3pt}\par
(3) For a simple generalized Witt algebra $W=W(A,T,\phi)$, the second
cohomology
group $H^2(W,\F)=0$ if ${\rm dim\ssc,}T\ge2$, and if $T=\F\partial$ is
1-dimensional, then $H^2(W,\F)$ is 1-dimensional, spanned by the
\nob{cohomology}
class $[\psi]$ where $\psi:W\times W\rar \F$ is the 2-cocycle defined by
$\psi(t^x\partial,t^y\partial)=\d_{x+y,0}\partial(x)^3,\,x,y\in A.$}
\qed\par
In this paper we are interested in the 1-dimensional universal central extension
of the simple algebras $W(A,T,\phi)$ in the case that $T=\F\partial$ is 1-dimensional.
In fact this leads us to the following definition.
\par\ni
{\bf Definition 2.2}.
Let $M$ be an abelian additive subgroup of $\F$, the {\it generalized
Virasoro algebra} $\Vir[M]$ is defined to be the Lie algebra with
$\F$-basis $\{L_\mu,c\,|\,\mu\in M\}$, subject to the following commutation
relations:
$$
\matrix{
[L_\mu,L_\nu]=(\nu-\mu)L_{\mu+\nu}+{1\over12}(\mu^3-\mu)
\d_{\mu+\nu,0}{\ssc\,}c,\vs{2pt}\hfill\cr
[c,L_\mu]=[c,c]=0,\hfill\cr
}
\eqno(2.3)$$
for $\mu,\nu\in M$.
\qed\par
It is straightforward to verify that $\Vir[M]$
is $M$-graded
$$
\matrix{
\Vir[M]=\oplus_{\mu\in M}\Vir[M]_\mu,\vs{3pt}\hfill\cr
\Vir[M]_\mu
=\left\{
\matrix{\F L_\mu,\ \,\,{\rm if }\,\,\mu\ne0,\vs{2pt}\hfill\cr
         \F L_0\oplus \F c,\ \,\,{\rm if }\,\,\mu=0,\vs{2pt}\hfill\cr
}\right.
\vs{3pt}\hfill\cr
[\Vir[M]_\mu,\Vir[M]_\nu] \subseteq\Vir[M]_{\mu+\nu}.\hfill\cr
}
\eqno(2.4)$$
{}From the above theorem we see that $\Vir[M]\simeq  \Vir[M']$ if  and only if
there exists $a\in\F^*$ such that $M'=aM$. It is clear that $\Vir[\Z]$ is the ordinary
Virasoro algebra, and if $M$ is of rank $n$, $\Vir[M]$ is a higher rank
Virasoro
algebra defined and studied in [PZ]. We can easily deduce the following theorem
from Theorem 2.1.
\par\ni
{\bf Theorem 2.3}. {\it (i) For  any $\chi\in{\rm Hom}(M,\F^*)$, the mapping
$$
\matrix{
\varphi_{\chi}\, :\,\Vir[M]\to\Vir[M],\vs{2pt}\hfill\cr
L_x\mapsto  \chi(x)L_x,\,\f x\in M; \,\,c\mapsto c\hfill\cr
}$$
is an  automorphism of  $\Vir[M]$.

\vs{-3pt}\par
(ii) For  any $a \in S(M):=\{\a\in\F\,|\,\a M=M\}$, the mapping
$$
\matrix{
\varphi'_{a}\, :\,\Vir[M]\to\Vir[M],\vs{2pt}\hfill\cr
L_x\mapsto  a^{-1}L_{ax},\,\f x\in M; \,\,c\mapsto a^{-1}c
\hfill\cr}
$$
is an  automorphism of  $\Vir[M]$.

\vs{-3pt}\par
(iii)  {\rm Aut}($\Vir[M])\simeq {\rm Hom}(M,\F^*)\times\hskip -5.3pt
\rb{1.8pt}{$\ssc|$}{\ssc\,}S(M)$.}\qed

\par
In  Theorem 3 of [PZ], mappings in  (a) are not automorphisms of  $\Vir[M]$,
mappings in (c) should be modified to (ii) of the above theorem. The reason
of those errors are the inaccuracy of Lemma 1 which we shall discuss in the next
section.
\par\
\par\
\par
\cl{\bf\S 3. Finite dimensional subalgebras of $\Vir[M]$}
\par
In this section we shall discuss finite dimensional subalgebras of $\Vir[M]$.
It suffices to consider this problem for only the centerless algebra
$ \tilde{\Vir}[M]=\Vir[M]/\F c$.
\par\ni
{\bf Lemma 3.1}. {\it Let $T$ be any finite dimensional subalgebras of
$ \tilde{ \Vir}[M]$. Then $\dim T\le3$. If $\dim T=3$, then there exists
a nonzero $n\in M$ such that $T=\F L_n+\F L_0+\F L_{-n}$.}

\par\ni
{\it Proof}. On the contrary, suppose $\dim T\ge4$, and
$X_i\,\,(i=1,2,3,4)$ are four linearly independent elements in $T$.
We choose a total ordering ``$\le$"on $M$ compatible with the group structure.
For an element $X=\sum_{\mu\in M}a_\mu L_\mu\in\tilde{\Vir}[M]$, we
define supp$(X)=\{\mu\in M\,|\,a_\mu\ne0\}$.
Let max\{supp$(X_i)$\}$=x_i$.
We may assume that $x_1>0$. If $x_2>0$ also, and $x_1\ne x_2$,
then $X_1,X_2$ generate an infinite dimensional subalgebra of $ \tilde{ \Vir}[M]$,
a contradiction. So $x_2=x_1$. By subtracting a suitable multiple of $X_1$
from $X_2$, we get $X'_2$ with  max\{supp$(X'_2)$\}$\le0$. In this way
we may assume that max\{supp$(X_i)$\}$\le0$ for $i=2,3,4$.
We consider min\{supp$(X_i)$\}$=x_i\le0$ for $i=2,3,4$.
Assume $x_2<0$. By subtracting suitable multiples of $X_2$
from $X_3, X_4$ respectively, we get two linearly independent
elements $X'_3,X'_4$ such that
$X'_3,X'_4\in\F L_0$, a contradiction. Thus $\dim T\le3$.

If $\dim T=3$, from the above discussion we can choose
a basis $\{Y_1,Y_2,Y_3\}$ such that supp($Y_1)\ge0$,
supp($Y_2)=0$, supp($Y_3)\le0$. Then there exists
a nonzero integer $n$ such that $T=\F L_n+\F L_0+\F L_{-n}$.\qed
\par\ni
{\bf Lemma 3.2}. {\it (i) Let $x\in M\setminus\{0\}$. Then for any
positive integer $n$ and any $\a\in\F$, the following two elements
span a two dimensional subalgebra of $  \tilde{\Vir}[M]$:
$$X=L_0+\a L_{-x},\,\,\, Y={\rm exp}(\a\,{\rm ad\ssc\,} L_{-x}) L_{nx}.$$

(ii) If $T$ is a two dimensional subalgebra of $  \tilde{ \Vir}[M]$, then     $T$ is
not abelian. If further $X=L_0+\a L_{-x}\in T$,
then there is an element $Y\in T$ such that
$Y={\rm exp}(\a\,{\rm ad\ssc\,} L_{-x}) L_{nx}$
for  some integer $n>0$.}

\par\ni
{\it Proof}. Part (i) is easy to verify. We omit the details.

(ii) It is clear that $T$ is not abelian. Since $\dim T=2$, there exists another
$Y\in  T\setminus \F X$ such that
$[X, Y]=aX+bY$ for some $a,b\in\F$. Further we can change $Y$, such that
$[X, Y]=bY\ne0$. Applying exp(ad$(-aL_{-x}))$ to  $[X, Y]=bY$, we get
$$[L_0, {\rm exp(ad}(-bL_{-x}))Y]= {\rm exp(ad}(-bL_{-x}))Y,$$
where exp(ad$(-bL_{-x}))Y$ can be an infinite sum. Thus
there exists a positive integer $n$ such that
exp(ad$(-bL_{-x}))Y=b'L_{nx}$, i.e., $Y=b'{\ssc\,}{\rm exp}(b\,{\rm ad\ssc\,} L_{-x}) L_{nx}$.\qed

\par
The above Lemma does not determine all 2-dimensional subalgebras
of $  \tilde{ \Vir}[M]$.
Let $X=(3/16)L_0+L_1+L_2,$  $Y=(3/16)L_0+(1/16)L_{-1}+(1/16^2)L_{-2}$, then
$\F X+\F Y$ is a subalgebra.  This subalgebra and the ones in Lemma 3.2
are not included in [PZ, Lemma 1] and [PZ, Theorem 4].
There are certainly some other subalgebras. It is not easy to determine all 2-dimensional  subalgebras of
$\tilde{ \Vir}[M]$.
Here we are not able to solve this problem,
but we can reduce this problem to
determining all 2-dimensional subalgebras of $  \tilde{ \Vir}[\Z]$.

\par\ni
{\bf Lemma 3.3}. {\it Let  $\{X,Y\}$ be a basis of a two dimensional subalgebra
$T$ of $  \tilde{ \Vir}[M]$. then
${\rm span\{supp}(X),{\rm supp}(Y)\}$ is a group of rank 1,
i.e., isomorphic to $\Z$.}

\par\ni
{\it Proof}. It is clear  that span$\{$supp$X$,  span$Y\}$ is a free group of finite rank.
We may assume that $[X, Y]=bY\ne0$. Choose a total ordering ``$\le$"on
$M$ compatible with the group structure.
If max\{supp$(Y)\}>0$, and max\{supp$(X)\}>0$, they  must be equal.
By subtracting a suitable multiple of $Y$
from $X$, we get $X'$ with  max\{supp$(X')\}=0$.
Then $[X', Y]=bY\ne0$. If supp$(X')=\{0\}$, we obtain that supp$(Y)$ is a singleton,
the lemma is true in this case. Suppose
min\{supp$(X')\}<0$. If  min\{supp$(Y)\}<0$,
they  must be equal since otherwise dim$T>2$.
By subtracting a suitable multiple of $X$ from $Y$,
we get $Y'$ with  min\{supp$(Y')\}=0$.
Since
${\rm span\{supp}(X),{\rm supp}(Y)\}={\rm span\{supp}(X'),{\rm supp}(Y')\}$,
if it is a group of rank $>1$, then we can choose another
ordering $\le'$ of $M$ such that one of the following is true:
(1) max\{supp$(Y')\}>0$ and  max\{supp$(X')\}>0$,  (2) min\{supp$(X')\}<0$,
and min\{supp$(Y')\}<0$.  In either cases we deduce that dim$T>2$,  a contradiction. Thus
our lemma follows.\qed
\par\
\par\
\par
\cl{\bf\S 4. Modules of the intermediate series}
\par
Let $\Vir[M]$ be the generalized Virasoro algebra defined in Definition
2.2. A {\it Harish-Chandra module}  over $\Vir[M]$ is a module $V$ such that
$$
V=\oplus_{\l\in\F} V_\l,\ V_\l=\{v\in V\,|\,L_0v=\l v\},
{\rm dim\ssc\,}V_\l<\infty,\,\f \l\in \F.
\eqno(4.1)$$
\par\ni
{\bf Definition 4.1}. {\it A module of the intermediate series}
over \,$\Vir[M]$ \,is an indecomposable \nob{Harish-Chandra}
module $V$ such that ${\rm dim\ssc\,}V_\l\le1$ for all $\l\in\F$.
\qed\par\ni
For any $a,b\in\F$, as in the Virasoro algebra case, one can define
the following
three series of $\Vir[M]$-modules $A_{a,b}(M),A_{a}(M),B_{a}(M)$, they all have
basis $\{v_\mu\,|\,\mu\in M\}$ with actions $c{\ssc\,}v_{\mu}=0$ and
$$
\matrix{
A_{a,b}(M):\hfill&
  L_\mu v_\nu=(a+\nu+\mu b)v_{\mu+\nu},\hfill&\vs{2pt}\hfill\cr
A_ {a}(M)  :\hfill&
  L_\mu v_\nu=(\nu+\mu)v_{\mu+\nu},\ \nu\ne0,\hfill&
  L_\mu v_0=\mu(\mu+a)v_\mu,\vs{2pt}\hfill\cr
B_  {a}(M)  :\hfill&
  L_\mu v_\nu=\nu{\ssc\,}v_{\mu+\nu},\ \nu\ne-\mu,\hfill&
  L_\mu v_{-\mu}=-\mu(\mu+a)v_0,\hfill\cr}
\eqno\matrix{\vs{2pt}(4.2{\rm a})\cr\vs{2pt}(4.2{\rm b})\cr(4.2{\rm c})\cr}
$$
for all $\mu,\nu\in M$. We use $A'_{a,b}(M)$, $A'_{a}(M)$,
$B'_{a}(M)$ to denote the nontrivial
subquotient of $A_{a,b}(M)$, $ A_{a}(M)$, $  B_{a}(M)$ respectively.
As results in [KR], $A'_{a,b}(M)\ne A_{a,b}(M)$
if and only if $a\in   M$ and $b=0$, or $a\in   M$ and $b=1$.
Note that we made a slight change in the modules
$A_ {a}(M)$, $B_ {a}(M)$.  But the notation here is obviously neater and simpler
than the old ones. You can see  this from the following theorem which is similar to [S, Prop.2.2].
\par\ni
{\bf Theorem 4.2}. {\it Among the $\Vir[M]$-modules  $A_{a,b}(M)$, $A_ {a}(M)$,
$B_  {a}(M) $ for $a,b\in\F$,
and their nontrivial subquotients, we have only the following module
isomorphisms:

(i) $A_{a,b}(M)\simeq  A_{a',b}(M)$ if $a-a'\in M$,

(ii) $A_{a,0}(M)\simeq  A_{a',1}(M)$ for $a\notin M$ with $a-a'\in M$,

(iii) $A'_{a,b}(M)\simeq  A'_{a',b}(M)$ if $a-a'\in M$,

(iv) $A'_{a,0}(M)\simeq  A'_{a',1}(M)$ for $a\in\F$ with $a-a'\in M$,

(v) $A'_{a}(M)\simeq  B'_{b}(M)\simeq A_{0,0}(M)$ for $a,b\in\F$.}\qed

\par
The following lemma is quite clear.
\par\ni
{\bf Lemma 4.3}. {\it Suppose $M_0\subseteq M$ is a subgroup of $M$,
and $V$ is a Harish-Chandra module over $\Vir[M]$ with a weight $\a$. If
$$V\simeq {\cases{ A_{a,b}(M), \,\,\, \a=a,\cr  B_  {a}(M), \,\,\a=0, \cr A_ {a}(M),\,\,\a=0, \cr A'_{0,0}(M)\oplus\F v_0, \,\,\a=0, \cr
A_{0,0}(M),\,\,\a=0, }}$$
for some $a,b\in\F$, then there exists $x_0\in M$ such that
$$V(\a+x_0,M_0):=\oplus_{z\in M_0}V_{\a+x_0+z}
\simeq {\cases{ A_{a,b}(M_0), \cr  B_  {a}(M_0 ), \cr A_ {a}(M_0 ),\cr A'_{0,0}(M_0 )\oplus\F v_0, \cr A_{0,0}(M_0 ),}}$$
respectively, and in the last four cases $\a+x_0\in M_0$.}\qed

\par
Su [S] proved the following
\par\ni
{\bf Theorem 4.4}. {\it Suppose $M\subseteq \F$ is an additive subgroup of rank $n$,
and $\F$ is algebraically closed of characteristic 0.
A module of the intermediate series over $\Vir[M]$
is isomorphic to one of
the  following: $A_{a,b}(M)$, $A_ {a}(M)$,
$B_{a}(M)$ for $a,b\in\F$, and their nonzero subquotients.}\qed

\par
It is natural to ask: are these three series of modules the only modules
of the intermediate series over the generalized Virasoro algebra $\Vir[M]$?
You will see the affirmative answer. Before we proceed the
proof, we need an auxiliary theorem. Similar to Theorem 2.3 of [OZ5]
(where the result is the same as in the following theorem but
only for $\Vir[\Z]$) we have

\par\ni
{\bf Theorem 4.5}. {\it Suppose $M\subseteq \F$ is an additive subgroup of rank $n$.
Let $V=\oplus_{j\in M}V_j$ be a $M$-graded $\Vir[M]$-module
with $\dim V>1$ and $\dim V_j\le1$
for all $j\in M$. Suppose there exists $a\in\F$ such that
$L_0$ acts on $V_j$ as the scalar $a+j$.
Then $V$ is isomorphic to one of the following for appropriate $\a,\b\in\F$:
(i) $A'_{\a,\b}(M)$, (ii) $A'_{0,0}\oplus\F v_0$ as direct sum of
$\Vir[M]$-modules, (iii) $A_{\a}(M)$, (iv) $B_{\a}(M)$, (v) $A_{0,0}(M)$.}

\par\ni
{\it Proof}. We may assume that $\Z\subset M$.
Denote the algebraically closed extension of $\F$ by $\ol{\F}$.
Write $\Sl_2=\ol{\F}d_1+\ol{\F}d_0+\ol{\F}d_{-1}$ which is the 3-dimensional
simple Lie algebra.
Consider the $\ol{\F}$ extensions $\ol V=\ol{\F}\otimes_{\F}V$,
$\ol{\Vir}[M]=\ol{\F}\otimes_{\F}\Vir[M]$, $\ol V_i=\ol{\F}\otimes_{\F}V_i$.
It is well known that any submodule of a $M$-graded module is still
$M$-graded.
\vskip 5pt
{\bf Claim 1} {\it $\ol V$ has no finite dimensional $\ol {\Vir}[M]$-submodule
of dimension $>1$.}
\vskip 5pt
Suppose, to the contrary, that $U$ is a finite dimensional
$\ol {\Vir} [M]$-submodule of $\ol V$ of \nob{dimension} $>1$.
{}From [S, Theorem 2.1], we see that each irreducible subquotient
of $U$ must be a \nob{1-dimensional} trivial $\ol{\Vir}[M]$-module.
Contrary to the fact that
 each weight space of $\ol V$ with respect to $L_0$ is 1-dimensional.
So Claim 1 holds.

\par
This claim implies that if $U$ is a subquotient of $\ol V$ then $\dim U=1$ or $\infty$.

\par
{\bf Claim 2} {\it If $U$ is a $\ol {\Vir} [M]$-subquotient of $\ol V$ and
$\dim U=\infty$, then $U\cap V_{a+x}\ne0$ for all $x\in M$ with $x+a\ne0$.}
\vskip 5pt
Suppose, to the contrary, that $U\cap V_{a+x_0}=0$ for $x_0\in M$ with
$x_0+a\ne0$. Choose $x_1\in M$ such that $U\cap V_{a+x_1}\ne0$.
Consider the $\ol{\Vir} [\Z(x_1-x_0)]$-module
$U'=\oplus_{k\in\Z}V_{a+x_0+k(x_1-x_0)}$.
By [OZ5, Theorem 2.3] we must have $U'\simeq A'_{0,0}(\Z(x_1-x_0))$. This
forces $a+x_0=0$, contrary to the choice of $x_0$. Thus
Claim 2 follows.

\par
We shall determine $V$ case by case.

{\it Case 1:} $\bvv$ is irreducible.

It follows from [S, Theorem 2.1] that $\bvv
\simeq A'_{\a,\b}(M)$ for appropriate $\a,\b\in
\ol {\F}$. Note that $\a\in {\F}$. If $\bvv\simeq A'_{0,0}(M)$, we can choose a basis
$\{v_i\,|\,i\in M\,\,{\rm with }\,\,
i\ne0\}$ for $\bvv$ such that $L_iv_j=jv_{i+j}$
with $v_0=0$. There exists $\g\in\ol {\F}$ so that $\g v_1\in V\setminus\{0\}$.
Then $\g v_i\in V$. Thus $V\simeq A'_{0,0}(M)$.
If $\bvv\not\simeq A'_{0,0}(M)$, we can choose a basis
$\{v_i\,|\,i\in M\,\}$ for $\bvv$ such that
$$L_iv_j=(j+\a+i\b)v_{i+j}.\eqno(4.3)$$
There exists an integer $i_0$ so that $L_{\pm1}v_i\ne0$ for all integer
$i\ge i_0$.
We may assume that $\g v_{i_0}\in V$ where $\g\in\ol{\F}\setminus\{0\}$.
Define
$$w_{i+i_0}=\g L_1^iv_{i_0}=\g\sum_{j=0}^{i-1}(\a+\b+j+i_0)v_{i+i_0},\f i\ge0.$$
Then $w_i\in V$ for all $i\ge i_0$, and
$L_jw_i\in V$ for all $j\ge-2, i> i_0+2$. Using (4.3) we see that
$$
\matrix{
L_{-1}w_{i+i_0}=(\a+\b+i_0+i-1)(\a-\b+i_0+i)w_{i+i_0-1},\vs{2pt}\hfill\cr
L_{-2}w_{i+i_0}=(\a+\b+i_0+i-1)(\a+\b+i_0+i-2)
(\a-2\b+i_0+i)w_{i+i_0-1}.\hfill\cr}
$$
It follows that
$$
\matrix{
(\a+\b+i_0+i-1)(\a-\b+i_0+i)\in {\F},\f i\ge2,\vs{2pt}\hfill\cr
(\a+\b+i_0+i-1)(\a+\b+i_0+i-2)
(\a-2\b+i_0+i)\in {\F},\f i\ge2.\hfill\cr}
$$
We can deduce that $\b(\b-1), -2\b(\b-1)(\b-2)\in {\F}$.
Thus we get $\b\in {\F}$. Therefore if we replace $v_i$
with $\g v_i$ for all $i\in\Z$, then $v_i\in V$ and
(4.3) still holds. Thus $V\simeq A'_{\a,\b}(M)$ with $\a,\b\in \F$.

\vskip 5pt
{\it Case 2:} $\bvv$ is reducible.

Let $U$ be a proper submodule of $\bvv$. From Claims 1 and 2 we know that
$U$ is irreducible and that $\dim U=1$ or $\infty$.

\vskip 5pt
{\it Subcase 1:} $\dim U=1$.

Setting $U=\bvv_{i_0}$, we see that $L_0\bvv_{i_0}=0$.
Then $\bvv/U$ is irreducible. From Case 1
we must have $\bvv/U\simeq \ol A'_{0,0}(M):=\ol{\F}\otimes_{\F}A'_{0,0}(M)$.
If $\bvv$ is decomposable, then
$\bvv=\ol A'_{0,0}(M)\oplus \bvv_0$ as $\ol{\Vir}$-modules, which implies that
$V=A'_{0,0}(M)\oplus  V_0$ as Vir-modules.
If $\bvv$ is indecomposable, by [S, Theorem 2.1] we have $\bvv\simeq
 \ol A_{0,0}(M)$, or $\bvv\simeq\ol B_{\a}(M)$ for
appropriate $\a\in \ol {\F}$. In the first case we have $V\simeq
 A_{0,0}(M)$,
In the second case, by changing the basis we can see that,
in fact,
$\a\in \F$ and $V\simeq B_{\a}(M)$.

\vskip 5pt
{\it Subcase 2:} $\dim U=\infty$.

\vskip 5pt
Claim 1 tells us that
$\dim(\bvv/U)=1$, and $\bvv/U=\bvv_0$. This implies that $L_0V_{0}=0$ and
$\bvv=U\oplus \ol V_{0}$ as subspaces. We deduce that $U$ is irreducible.
By Case 1 we have $U\simeq  A'_{0,0}(M)$. If $\bvv$ is decomposable,
then $\bvv=\ol  A'_{0,0}(M)\oplus\bvv_{i_0}$ as submodules. And
it follows that $V= A'_{0,0}(M)\oplus V_{0} $  as Vir-submodules.
If $\bvv$ is indecomposable,
then [S, Theorem 2.1] ensures  $\bvv\simeq
 \ol A_{0,1}(M)$, or $\bvv\simeq A_{\a}(M)$ for
appropriate $\a\in \ol{\F}$. In the first case we have $V\simeq
 A_{0,1}(M)$,
In the second case, by changing the basis we can see that,
in fact,
$\a\in {\F}$ and $V\simeq A_{\a}(M)$.\qed
\par
Now we are ready to classify  modules of the intermediate series over
$\Vir[M]$ for any $M$ and any field $\F$ of characteristic 0.

\par\ni
{\bf Theorem 4.6}. {\it A module of the intermediate series over
$\Vir[M]$ is isomorphic to one of
the  following: $A_{a,b}(M)$, $A_ {a}(M)$,
$B_{a}(M)$ for $a,b\in\F$, and their nonzero subquotients.}

\par\ni
{\bf Proof.}
Suppose $V$ is a module of the intermediate series. As it is
indecomposable,
there exists some $a\in\F$ such that
$$
V=\oplus_{\mu\in M}V_{\mu+a}.
\eqno(4.4)$$
If $a\in M$, we choose $a=0$.

For $x\in M$ and any subgroup $N$ of $M$, we define
$$V(a+x,N)=\oplus_{z\in N}V_{a+x+z}.\eqno(4.5)$$
It clear that $V(a+x,N)$ is a $\Vir[N]$-module.

\par
If $\dim V=1$, then $V$ is the  1-dimensional trivial $\Vir[M]$-module.
Next we shall always assume that $\dim V>1$.
\par
If $M\simeq \Z$, [OZ5, Theorem 2.3] ensures that  Theorem 4.6 holds.
Next we shall always assume that $M\not\simeq \Z
$.
\par
{\bf Claim 1}. {\it $V_{a+x}\ne0$ for all $x\in M$ with $a+x\ne0$.}

\par
Suppose, to the contrary, that $V_{a+x_0}=0$ for $x_0\in M$ with
$x_0+a\ne0$. Choose $x_1\in M$ such that $V_{a+x_1}\ne0$.
Consider the $\Vir[\Z(x_1-x_0)]$-module $U=\oplus_{k\in\Z}V_{a+x_0+k(x_1-x_0)}$.
By [OZ5, Theorem 2.3] we must have $U\simeq A'_{0,0}(\Z(x_1-x_0))$. This
forces $a+x_0=0$, contrary to the choice of $x_0$. Thus
Claim 1 follows.

\par
{}From Claim 1 and Theorem 4.5 we know that $c$ acts trivially on $V$.

\par
Now we fix an arbitrary nonzero $z_0\in M$, let $\{x_i\,|\,i\in I\}$ be the set of all
representatives of cosets of $\Z z_0$ in $M$. If $a\in M$, we choose $x_0=0$.
It is clear that $V=\oplus_{i\in I} V(a+x_i,\Z z_0)$, see (4.5) for the definition.
{}From Claim 1 and Theorem 4.5 we know that, for any $i\in I\setminus\{0\}$, there exists $b_i\in\F$ such that
$$V(a+x_i, \Z z_0 )\simeq  A_{a+x_i,b_i}(\Z z_0 ),\eqno(4.6) $$
and
$$V(a+x_0, \Z z_0 )\simeq {\cases{ A_{a+x_0,b_0}(\Z z_0 ),\cr  B_  {\a}(\Z z_0 ), \cr A_ {\a}(\Z z_0 ), \cr A'_{0,0}(\Z z_0 )\oplus\F v_0, \cr
A'_{0,0}(\Z z_0 ),}}$$
for some $b_0, \a\in\F$.

\par
Suppose $M_0$ is a maximal subgroup of $M$ satisfying the  following
two conditions:

(C1) Let $\{x_i\,|\,i\in I'\}$ be the set of all
representatives of cosets of $M_0$ in $M$ with $0\in  I'$ and $x_0=0$.
Then for all $i\in  I'\setminus\{0\}$,
$$V(a+x_i, M_0 )\simeq   A_{a+x_i,b_i}(M_0 ),$$
for some $b_i\in\F$;

(C2) $$V(a+x_0, M_0 )\simeq {\cases{ A_{a+x_0,b_0}(M_0 ),\cr  B_  {\a}(M_0  ), \cr A_ {\a}(M_0 ), \cr A'_{0,0}(M_0 )\oplus
\F v_0, \cr
A'_{0,0}(M_0 ),}}$$
for some $b_0, \a\in\F$.

\par
{}From the above discussion we know that such an $M_0$ exists.
It suffices to show  that $M_0=M$.
\par
Otherwise we suppose $M_0\ne M$. So $|I'|>1$, and $x_1\in  M\setminus M_0$ for $1\in I'$.
Denote $M_1=M_0+\Z x_1$.
Let $\{y_i\,|\,i\in J\}$ be the set of all
representatives of cosets of $M_1$ in $M$  with $0\in  J$ and $y_0=0$,
and  $\{ix_1\,|\,i\in K\}$ be the set of all
representatives of cosets of $M_0$ in $M_1$, where $K\subset\Z$.

\par
{\bf Claim 2}. {\it For any fixed  $j\in J$ we have
$$V(a+y_j, M_1 )\simeq {\cases{ A_{a+y_j,\a}(M_1),\cr  B_  {\a}(M_1), \cr A_ {\a}(M_1 ), \cr A'_{0,0}(M_1 )\oplus
\F v_0, \cr
A'_{0,0}(M_1),}}$$
for some $\a\in\F$.}

\par
We shall show this claim in two cases.

\par
{\it Case 1}: For all $i\in K$,
$$V(a+ix_1+y_j, M_0 )\simeq   A_{a+ix_1+y_j,b_i}(M_0 ),\eqno(4.7)$$
for some $b_i\in\F$.

\par
{}From Theorem 4.2, we have two choices for  some $b_i$.
{}From (4.7) we deduce that
$$V(a+ix_1+y_j,\Z z_0)\simeq A_{a+y_j+ix_1,b_i}(\Z z_0).\eqno(4.8)$$
For any $z_0\in M_0\setminus\{0\}$, rank$(\Z z_0+\Z x_i)=1$
or $2$. Then for all $i\in K$,
[S, Theorem 2.1] and Lemma 4.3 ensure that
$V(a+ix_1+y_j,\Z z_0+\Z x_1)=V(a+y_j,\Z z_0+\Z x_1)
=\oplus_{z\in\sZ z_0+\sZ x_1}V_{a+y_j+z}\simeq A_{a+y_j,b}(\Z z_0+\Z x_1)$ for
some $b\in\F$. Thus for all $i\in K$,
$$
\matrix{
V(a+ix_1+y_j,\Z z_0)\simeq A_{a+y_j+ix_1,b}(\Z z_0),\vs{3pt}\hfill\cr
V(a+y_j,\Z x_1)\simeq A_{a+y_j,b}(\Z x_1),\hfill\cr
}
\eqno\matrix{\vs{3pt}(4.9)\cr(4.10)\cr}$$
Combining  (4.8) with (4.9), we can choose $b_i$ such that
  $b=b_i$ for all $i\in \Z$, i.e.,
$$V(a+ix_1+y_j, M_0 )\simeq   A_{a+ix_1+y_j,b}(M_0 ),\eqno(4.11)$$

Using (4.10) and (4.11), we choose a basis $v_{y_j+z}\in V_{a+y_j+z}$
for all $z\in M_1$ and  all $i\in K$ as follows:

(N1) Choose $v_{kx_1+y_j}\in V_{a+kx_1+y_j}$
for all $k\in \Z$ such that
$$L_{ix_1}v_{y_j+kx_1}=(a+y_j+kx_1+ix_1b)v_{y_j+kx_1+ix_1}, \f i,k\in\Z;$$

(N2) For  any  $i\in K$, choose $v_{ix_1+y_j+z}\in V_{a+ix_1+y_j+z}$
for all $z\in M_0\setminus\{0\}$ such that
$$L_{y}v_{y_j+ix_1+z}=(a+y_j+ix_1+z+by)v_{y_j+ix_1+z+y},\f y,z\in M_0.$$

\par
Applying  Lemma 4.3 to (4.10) gives
$$V(a+y_j, \Z z_0+\Z x_1)\simeq   A_{a+y_j,b}(\Z z_0+\Z x_1).\eqno(4.12)$$
Then we can choose a basis $v'_{y_j+z}\in V_{a+y_j+z}$
for all $z\in \Z z_0+\Z x_1$ with $v'_{y_j}= v_{y_j}$ such that
$$L_{y}v'_{y_j+ix_1+z}=(a+y_j+ix_1+z+by)v'_{y_j+ix_1+z+y}.$$
Combining this with (N1) and (N2), we deduce that
$v'_{y_j+z}=v_{y_j+z}$ for all $z\in \Z z_0+\Z x_1$, i.e.,
$$L_{y}v'_{y_j+z}=(a+y_j+z+by)v'_{y_j+z+y},\,\,\f y,z\in \Z z_0+\Z x_1.$$
Since $z_0\in M_0$ is arbitrary, the above equation holds for all $y,z\in M_1$.
This yields that $V(a+y_j,M_1)\simeq A_{a+y_j,b}(M_1)$. So Claim 2 is true for this case.

\par
{\it Case 2}: $a=0$ and
$$V(0, M_0 )=\oplus_{z\in M_0}V_z
\simeq {\cases{ B_  {\a}(M_0), \cr A_ {\a}(M_0 ), \cr A'_{0,0}(M_0 )\oplus\F
v_0,\,\,\,{\rm or } \cr A'_{0,0}(M_0 ),}}\eqno(4.13)$$
for some $\a\in\F$.

\par
In this case, $j=0,\,\,y_j=0$.
It is clear that, for any  $i\in K\setminus\{0\}$, since $ix_1\notin M_0$
we have
$$V(ix_1,M_0)=A_{ix_1,b_i}(M_0 ).\eqno(4.14)$$
For any $z\in M_0\setminus\Z x_1$, rank$(\Z z+\Z x_i)=1$ or $2$.
If $M_0\subset \Z x_1$, then $M_1\simeq \Z$, Claim 2 is automatically true.
Now suppose $z_0\in M_0\setminus\Z x_1$ is a fixed arbitrary element.
Then for all $i\in K$,
[S, Theorem 2.1], Lemma 4.3 and (4.13) ensure that
$$V(ix_1,\Z z_0+\Z x_1)=V(0,\Z z_0+\Z x_1)\simeq
{\cases{ B_  {\a}(\Z z_0+\Z x_1), \cr A_ {\a}(\Z z_0+\Z x_1), \cr A'_{0,0}(\Z z_0+\Z x_1)\oplus
\F v_0, \,\,\,{\rm or } \cr
A'_{0,0}(\Z z_0+\Z x_1),}}$$
respectively. Thus
$$V(ix_1,\Z z_0)\simeq {\cases{ A_{ix_1,0}(\Z z_0) ,\cr A_{ix_1,1}(\Z z_0) , \cr A'_{ix_1,0}(\Z z_0),\,\,\,{\rm or } \cr
A'_{ix_1,0}(\Z z_0) ,}}\,\,\,\,\,\,\,\,\f i\in K\setminus\{0\},\eqno(4.15)$$
$$V(0,\Z z_0 )\simeq {\cases{ B_{\a}(\Z z_0 ), \cr A_{\a}(\Z  z_0 ), \cr A'_{0,0}(\Z  z_0 )\oplus
\F v_0, \,\,\,{\rm or } \cr
A'_{0,0}(\Z z_0 ),}}\eqno(4.16)$$
$$V(z_0,\Z x_1 )\simeq {\cases{ A_{z_0,0}(\Z x_1) ,\cr A_{z_0,1}(\Z x_1 ) , \cr A'_{z_0,0}(\Z z_0),\,\,\,{\rm or } \cr
A'_{z_0,0}(\Z x_1 ), }}\eqno(4.17)$$
respectively.
{}From (4.14) and (4.15) we can  choose $b_i$ ($i\in K\setminus\{0\}$) such that
$$b_i=b=0, 1, 0 ,\,\,\,{\rm or}\,\, 0\eqno(4.18)$$
corresponding to the four cases respectively.
Using (4.15), (4.16), (4.17), (4.18) and (4.14), we choose
a basis $v_{ix_1+z}\in V_{ix_1+z}$
for all $z\in M_0$ and  all $i\in K$ as follows:

(N1') Choose $v_{z_0+kx_1}\in V_{z_0+kx_1}$
for all $k\in \Z$ such that
$${\cases{L_{ix_1}v_{z_0+jx_1}=(z_0+jx_1 )v_{z_0+(i+j)x_1}, \,\,\f \,\,i,j\in\Z,\vs{2pt}\cr
L_{ix_1}v_{z_0+jx_1}=( z_0+ (i+j)x_1 )v_{z_0+(i+j)x_1},\,\,\f \,\, i,j\in\Z,\vs{2pt}\cr
L_{ix_1}v_{z_0+jx_1}=(z_0+ jx_1 )v_{z_0+(i+j)x_1},\,\,\f \,\,i,j\in\Z,\vs{2pt}\cr
L_{ix_1}v_{z_0+jx_1}=(z_0+ jx_1) v_{z_0+(i+j)x_1},\,\,\f\,\, i,j\in\Z,}}$$
respectively;

(N2') For  any  $j\in K$, choose $v_{jx_1+z}\in V_{jx_1+z}$
for all $z\in M_0\setminus\{z_0\}$ such that, if $j\in K\setminus\{0\}$,
$${\cases{L_{z}v_{z_0+jx_1+y }=(z_0+jx_1+y  )v_{z_0+jx_1+z+y }, \,\,\f \,\,y,z\in M_0,\vs{2pt}\cr
L_{z}v_{z_0+jx_1+y }=( z_0+ z+jx_1+y  )v_{z_0+jx_1+z+y },\,\,\f \,\, y,z\in M_0,\vs{2pt}\cr
L_{z}v_{z_0+jx_1+y }=(z_0+ jx_1+y  )v_{z_0+jx_1+z+y },\,\,\f \,\,y,z\in M_0,\vs{2pt}\cr
L_{z}v_{z_0+jx_1+y }=(z_0+ jx_1+y ) v_{z_0+jx_1+z+y },\,\,\f\,\, y,z\in M_0,}}$$
respectively; if $j=0$, $\f y,z\in M_0$,
$${\cases{L_{z}v_{y }=y v_{z+y }\, \,{\rm for}\,\,y+z\ne0,\,\,
L_zv_{-z}=z(z+\a)v_0,\vs{2pt}\cr
L_{z}v_{y }=( z+y  )v_{z+y }\,\,{\rm for}\,\,y\ne0,\,\,L_zv_{0}=z(z+\a)v_z,\vs{2pt}\cr
L_{z}v_{y }=yv_{z+y } \,\,{\rm for}\,\,y(y+z)\ne0,\vs{2pt}\cr
L_{z}v_{y }=y  v_{z+y }\,\,{\rm for}\,\,y(y+z)\ne0,}}$$
respectively.

Note that this choice of basis can always be done although we have sometimes
$M_0\cap \Z x_1\ne0$.
\par
For  any $z_1\in  M_0$,  (4.13) yields (4.13) with $M_0$ replaced into $\Z z_1$.
Applying  Lemma 4.3 gives
(4.13) with $M_0$ replaced into $\Z z_1+\Z x_1$.
Then we can choose a basis $v'_{z_0+z}\in V_{z_0+z}$
for all $z\in M_2:=\Z z_1+\Z x_1$ with $v'_{z_0}= v_{z_0}$ such that,
for all $y,z\in M_2$,
$${\cases{L_{z}v'_{z_0+y }=(z_0+y  )v'_{z_0+z+y } \,\,\,{\rm if}\,\,\,
z_0+y+z\ne0,\,\, L_zv'_{-z}=z(z+\a)v'_0, \vs{2pt}\cr
L_{z}v'_{z_0+y }=( z_0+ z+y  )v'_{z_0+z+y }\,\,\,{\rm if}\,\,\,
z_0+y\ne0,\,\, L_zv'_{0}=z(z+\a)v'_z,\vs{2pt}\cr
L_{z}v'_{z_0+y }=(z_0+ y  )v'_{z_0+z+y }, \vs{2pt}\cr
L_{z}v'_{z_0+y }=(z_0+y ) v'_{z_0+z+y },}}$$
respectively.
Combining this with (N1') and (N2'), we deduce that
$v'_{z_0+z}=v_{z_0+z}$ for all $z\in \Z z_1+\Z x_1$.
Since $z_1\in M_1$ is arbitrary, we  can deduce (4.13) with $M_0$ replaced by
$M_1$. So Claim 2 is  true for this subcase.
\par
Claim 2 is contrary to the choice of $M_0$.  Therefore we must have $M_0=M$.
This completes the proof of this theorem.
\qed
\par\
\par\
\par
\cl{\bf\S 5. Generalized super-Virasoro algebras}
\def\SVir{\mbox{\bf\rm SVir}}
\par
In this section, we shall first introduce the notion of the generalized
super-Virasoro algebras which generalizes the notion of the high rank
super-Virasoro algebras introduced in [S], and then determine
the modules of the intermediate series over the
generalized super-Virasoro algebras.
\par
Roughly speaking, a generalized super-Virasoro algebra is a Lie superalgebra
which is a nontrivial $\Z/2\Z$-graded extension of a generalized
Virasoro algebra.
Thus, suppose $\SVir[M]=\SVir_{\ol0}[M]\oplus\SVir_{\ol1}[M]$ is a
$\Z/2\Z$-graded extension of a generalized Virasoro algebra
$\Vir[M]$ such that $\SVir_{\ol0}[M]=\Vir[M]$ and $\SVir_{\ol1}[M]$ is
a nontrivial irreducible $\Vir[M]$-module of the intermediate series.
By Theorem 4.2 and Theorem 4.6, $\SVir_{\ol1}[M]$ is a subquotient module
of $A_{\a,b}(M)$ for some $\a,b\in\F$, so by rewriting (4.2a), there is a
subset $M'$ of $M$, where $M'=M$ if $\a\notin M,\,b\ne0,1$ or $M'=M\bs\{0\}$
if $\a=0,\,b=0,1$, such that there exists a basis $\{G_\mu\,|\,\mu\in\a+M'\}$
and
$$
\matrix{
[L_\mu,G_\nu]=(\nu+\mu b)G_{\mu+\nu},\vs{2pt}\hfill\cr
[c,G_\nu]=0,\hfill\cr}
\ \ \ \forall\ \mu\in M,\ \nu\in \a+M'.
\eqno(5.1)$$
Since we are considering a nontrivial extension, we have
$0\ne [G_\mu,G_\nu]\in\Vir[M]$ for some $\mu,\nu\in M'$.
As it has the weight $\mu+\nu$, thus we have $\mu+\nu\in M$ and so
$$
2\a\in M.
\eqno(5.2)$$
In general, since
$[\SVir_{\ol1}[M],\SVir_{\ol1}[M]]\subset\SVir_{\ol0}[M]=\Vir[M]$, we can write
$$
[G_\mu,G_\nu]=x_{\mu,\nu}L_{\mu+\nu}+\d_{\mu+\nu,0}y_{\mu}c,
\ \ \forall\ \mu,\nu\in\a+M',
\eqno(5.3)$$
for some $x_{\mu,\nu},y_{\mu}\in\F$.
Applying ${\rm ad\ssc\,}L_\l,\l\in M$ to (5.3), using (5.1) and definition
(2.3), we obtain that
$$
\matrix{
(\mu+\l b)x_{\mu+\l,\nu}+(\nu+\l b)x_{\mu,\nu+\l}=(\mu+\nu-\l)x_{\mu,\nu},
\vs{2pt}\hfill\cr
\d_{\mu+\l+\nu,0}(\mu+\l b)y_{\mu+\l}+\d_{\mu+\nu+\l,0}(\nu+\l b)y_\mu=
{1\over12}(\l^3-\l)\d_{\l+\mu+\nu,0}x_{\mu,\nu},
\hfill\cr
}
\eqno\matrix{\vs{2pt}(5.4{\rm a})\cr(5.4{\rm b})\cr}$$
holds for all $\l\in M,\,\mu,\nu\in\a+M'$.
\par
We shall consider (5.4a-b) in two cases.
\par
{\it Case 1:} Suppose $[\SVir_{\ol1}[M],\SVir_{\ol1}[M]]\subseteq\F c$, i.e.,
$x_{\mu,\nu}=0$ for all $\mu,\nu\in\a+M'$.
\par
By taking $\nu=\mu\in\a+M',\,\l=-2\mu\in M$ in (5.4b), we obtain
$$
\mu(1-2b)(y_{-\mu}+y_\mu)=0,\ \ \forall\ \mu\in \a+M'.
\eqno(5.5)$$
First suppose $b\ne{1\over2}$. Then (5.5) gives $y_{-\mu}=-y_\mu$.
On the other hand, from the definition of Lie superalgebras we have
$$y_{-\mu}=y_\mu, \ \ \forall\,\mu\in \a+M'.
\eqno(5.6)$$
This forces $y_\mu=0,\,\forall\,\mu\in \a+M'$ and so,
we obtain the trivial extension. Therefore $b={1\over2}$, and so $M'=M$.
Then in (5.4b), by taking $\l=-(\nu+\mu)\in M$, we obtain
$$
\mbox{$1\over2$}(\mu-\nu)(y_{-\nu}-y_\mu)=0,\ \ \forall\ \mu,\nu\in \a+M.
\eqno(5.7)$$
In particular, by setting $\nu=\a$, we have
$y_\mu=y_{-\a}$ for all $\a\ne\mu\in \a+M$ and by (5.6), we have
$y_\a=y_{-\a}$. This shows that $y_\mu=y\in\F$ must be a nonzero scalar.
By rescaling the basis $\{G_\mu\,|\,\mu\in\a+M\}$ if necessary, we can
suppose $y_\mu=1$ for all $\mu\in \a+M$ and so from (5.1) and (5.3), we have
$$
\matrix{
[L_\mu,G_\nu]=(\nu+{1\over2}\mu)G_{\mu+\nu},\ [c,G_\nu]=0,
\vs{2pt}\hfill\cr
[G_\nu,G_\l]=\d_{\nu+\l,0}c,\hfill\cr
}\ \ \ \forall\ \mu\in M,\,\nu,\l\in \a+M.
\eqno(5.8)$$
\par
It is immediate to check that
$\wt{\rm S}{\rm Vir}[M,\a]={\rm span}\{L_\mu,G_\nu,c\,|
\,\mu\in M,\nu\in \a+M\}$ with commutation relations (2.3) and (5.8)
defines a Lie superalgebra for any subgroup $M$ of $\F$ and $\a\in\F$ such
that $2\a\in M$.
\par
{\it Case 2:} Suppose $[\SVir_{\ol1}[M],\SVir_{\ol1}[M]]\not\subseteq\F c$, i.e.,
$x_{\mu,\nu}\ne0$ for some $\mu,\nu\in \a+M'$.
\par
Using the super-Jacobian identity,
$$
[G_\l,[G_\mu,G_\nu]]=[[G_\l,G_\mu],G_\nu]-[G_\mu,[G_\l,G_\nu]],
\ \ \forall\,\mu,\nu,\l\in \a+M',
\eqno(5.9)$$
by applying $G_\l$ to (5.3) and using (5.1), we have
$$
-(\l+(\mu+\nu)b)x_{\mu,\nu}
=(\nu+(\l+\mu)b)x_{\l,\mu}+(\mu+(\l+\nu)b)x_{\l,\nu},
\ \ \forall\,\l,\mu,\nu\in \a+M'.
\eqno(5.10)$$
Setting $\l=\nu=\mu\in \a+M'$ in (5.10), we obtain
$$
3\mu(1+2b)x_{\mu,\mu}=0,\ \ \ \forall\ \mu\in \a+M'.
\eqno(5.11)$$
If $b\ne-{1\over2}$, then
$$
x_{\mu,\mu}=0, \ \ \ \forall\ 0\ne\mu\in \a+M'.
\eqno(5.12)$$
By setting $\nu=\mu$ in (5.4a), and then substituting $\mu+\l$ by
$\nu$, using (5.12) and noting that $x_{\mu,\nu}=x_{\nu,\mu}$, we obtain
$$
2(\mu(1-b)+\nu b)x_{\mu,\nu}=0,\ \ \forall\ \mu,\nu\in \a+M'.
\eqno(5.13)$$
Suppose $x_{\mu_0,\nu_0}\ne0$ for some $\mu_0,\nu_0\in \a+M'$, then (5.13) and
$x_{\mu_0,\nu_0}=x_{\nu_0,\mu_0}$ give
$\mu_0(1-b)+\nu_0b=0=\nu_0(1-b)+\mu_0b$. From this, we have $\mu_0=-\nu_0$, thus
we obtain
$$
x_{\mu,\nu}=0,\ \ \forall\,\mu,\nu\in \a+M',\ \mu\ne-\nu.
\eqno(5.14)$$
But when $\mu=-\nu$, by taking $\l\ne0$ in (5.4a), we again obtain
$x_{\mu,-\mu}=0$. Thus $x_{\mu,\nu}=0$ for all $\mu,\nu\in \a+M'$,
a contradiction.
Therefore we obtain $b=-{1\over2}$. So, $M'=M$.
Then in (5.10), by letting
$\l=\mu$, we obtain $x_{\mu,\nu}=x_{\mu,\mu}$ for all $\mu,\nu\in \a+M$.
Therefore we have $x_{\mu,\nu}=x_{\nu,\nu}=x_{\a,\nu}=x_{\a,\a}$ is a nonzero
scalar. By rescaling basis $\{G_\mu\,|\,\mu\in \a+M\}$ if necessary, we can
suppose $x_{\mu,\nu}=2$ for all $\mu,\nu\in \a+M$.
\par
By letting $\nu=\mu\in \a+M$, $\l=-2\mu\in M$ in (5.4b), we obtain
$$
\mu y_\mu=-\mbox{$1\over12$}\mu(4\mu^2-1), \ \ \forall\,\mu\in \a+M.
\eqno(5.15)$$
Thus if $\a\notin M$, we obtain that
$$
y_\mu=-\mbox{$1\over3$}(\mu^2-\mbox{$1\over4$}),
\eqno(5.16)$$
holds for all $\mu\in \a+M$. If $\a\in M$, then (5.16) holds for all
$\mu\ne0$. To prove (5.16) also holds for $\mu=0$ when $\a\in M$,
setting $\nu=0,\l=-\mu\in M$ in (5.4b), we see this immediately.
Thus (5.1) and (5.3) have the following forms:
$$
\matrix{
[L_\mu,G_\nu]=(\nu-{1\over2}\mu)G_{\mu+\nu},\ [c,G_\nu]=0,
\vs{2pt}\hfill\cr
[G_\nu,G_\l]=2L_{\nu+\l}-{1\over3}(\nu^2-{1\over4})\d_{\nu+\l,0}c,\hfill\cr
}\ \ \ \ \forall\ \mu\in M,\,\nu,\l\in \a+M.
\eqno(5.17)$$
\par
It is immediate to check that
$\SVir[M,\a]={\rm span}\{L_\mu,G_\nu,c\,|
\,\mu\in M,\nu\in \a+M\}$ with commutation relations (2.3) and (5.17)
defines a Lie superalgebra for any subgroup $M$ of $\F$ and $\a\in\F$ such
that $2\a\in M$.
\par
Thus we have in fact proved
\par\ni
{\bf Lemma 5.1}
{\it Suppose $W=W_{\bar 0}\oplus W_{\bar 1}$ is a  Lie superalgebra such that
$W_{\bar 0}\simeq \Vir[M]$ and $W_{\bar 1}$ is an irreducible $\Vir[M]$-module of
the intermediate series, where $M$ is an additive subgroup of  $\F$.  Then $W$ is
$\wt{\rm S}{\rm Vir}[M,\a]$ or $\SVir[M,\a]$ for a suitable
$\a\in\F$ with $2\a\in M$.}
\qed\par
This result leads us the following definition
\par\ni
{\bf Definition 5.2}.
A {\it generalized super-Virasoro algebra} is a Lie superalgebra
$\wt{\rm S}{\rm Vir}[M,\a]$ defined by (2.3) and (5.8) or
$\SVir[M,\a]$ defined by (2.3) and (5.17), where $M$ is a subgroup of
$\F$, $\a\in\F$ with $2\a\in M$.
\qed\par
Now we are in a position to consider modules of the intermediate series
over the generalized super-Virasoro algebras.
Since $\wt{\rm S}{\rm Vir}[M,\a]$ is just a trivial extension of $\Vir[M]$
modulo the center $\F c$, the modules of the intermediate series over
$\wt{\rm S}{\rm Vir}[M,\a]$ are simply those modules over the generalized
Virasoro algebra $\Vir[M]$. Thus we shall be only interested in considering
$\SVir[M,\a]$.
\par
First we give the precise definition.
\par\ni
{\bf Definition 5.3}.
A $\Z/2\Z$-graded vector space $V=V_{\ol0}\oplus V_{\ol1}$ is called
a {\it module of the \nob{\it intermediate} series} over $\SVir[M,\a]$ if it
is an indecomposable $\SVir[M,\a]$-module such that
$$
\matrix{
\SVir_\si[M,\a]V_\tau\subset V_{\si+\tau},\,\f\si,\tau\in\Z/2\Z,\vs{2pt}\hfill\cr
{\rm dim}_{\F}(V_\si)_\l\le1,
\mbox{ where } (V_\si)_\l=\{v\in V_\si\,|\,L_0v=\l v\},\,\f\l\in\F,\,
\si\in\Z/2\Z.
\hfill\cr
}
\eqno\matrix{\cr(5.18)\cr\hfill\qed\cr}$$
\vs{-4pt}\par\def\SA{\mbox{$S\!A$}}\def\SB{\mbox{$S\!B$}}
As in [S], [S1], there exist three series of modules
$\SA_{a,b}(M,\a),$ $\SA_a(M,\a),$ $\SB_a(M,\a)$ for
$a,b\!\in\!\F$ comparable with those of the generalized
Virasoro algebra defined in (4.2a-c). \nob{Precisely},
$\SA_{a,b}(M,\a),\SA_a(M,\a)$ have basis
$\{v_\mu\,|\,\mu\in M\}\cup\{w_\nu\,|\,\nu\in\a+M\}$ and $\SB_a[M,\a]$
has basis $\{v_\nu\,|\,\nu\in \a+M\}\cup\{w_\mu\,|\,\mu\in M\}$ such that
the central element $c$ acts trivially and
\par\ni\hs{3pt}$
\matrix{
\SA_{a,b}(M,\a)\!:\!\!\hfill&\!\!
  L_\l v_\mu=(a+\mu+\l b)v_{\l+\mu},\hfill&\!\!
  L_\l w_\nu=(a+\nu+\l(b-{1\over2}))w_{\l+\nu},
 \vs{2pt}\hfill\cr
 &\!\!
  G_\eta v_\mu=w_{\eta+\mu},\hfill&\!\!
  G_\eta w_\nu=(a+\nu+2\eta(b-{1\over2}))v_{\eta+\nu},
 \vs{4pt}\hfill\cr
}$\hfill(5.19a)\par\ni\hs{3pt}$
\matrix{
\SA_ {a}(M,\a)\!:\!\!\hfill&\!\!
  L_\l v_\mu=(\mu+\l)v_{\l+\mu},\mu\ne0,\hfill&\!\!
  L_\l v_0=\l(\l+a)v_\l,\hfill&\!\!
   L_\l w_\nu=(\nu+{\l\over2})w_{\l+\nu},\hfill
 \vs{2pt}\cr
 &\!\!
  G_\eta v_\mu=w_{\eta+\mu},\mu\ne0,\hfill&\!\!
  G_\eta v_0=(2\eta+a)w_\eta,
&\!\! G_\eta w_\nu=(\nu+\eta)v_{\eta+\nu},\hfill\cr
}\!\!\!\!\!$\hfill(5.19b)\par\ni\hs{3pt}$
\matrix{
\SB_{a}(M,\a)\!:\!\!\hfill&\!\!
  L_\l v_\nu=(\nu+{\l\over2})v_{\l+\mu},\hfill&\!\!
  L_\l w_\mu=\mu w_{\l+\mu},\mu\ne-\l,\hfill&\!\!
  L_\l w_{-\l}=-\mu(\mu+a)w_0,
 \vs{2pt}\hfill\cr
 &\!\!
  G_\eta v_\nu=w_{\eta+\nu},\nu\ne-\l,\hfill&\!\!
  G_\eta v_{-\eta}=(2\eta+a)w_0,\hfill
&\!\! G_\eta w_\mu=\mu v_{\eta+\mu},\hfill&\!\!\cr
}\!\!\!\!$\hfill(5.19c)\par\ni
%
for all $\l,\mu\in M,\,\nu,\eta\in\a+M$.
\par
Now we are ready to prove following theorem.
\par\ni
{\bf Theorem 5.4}.
{\it A module of the intermediate series over $\SVir[M,\a]$ is
isomorphic to one of the following: $\SA_{a,b}(M,\a)$,
$\SA_{a}(M,\a)$, $\SB_{a}(M,\a)$ for $a,b\in\F$, and their
nonzero subquotients.}
\par\ni
{\bf Proof.}
Let $V=V_{\ol0}\oplus V_{\ol1}$ be a module of the intermediate series
over $\SVir[M,\a]$.
Since $V$ is \nob{indecomposable,} there exists some $a\in\F$
such that $V_{\ol0}=\oplus_{\l\in M}(V_{\ol0})_{a+\l}$
and $V_{\ol1}=\oplus_{\l\in M}(V_{\ol1})_{a+\a+\l}$.
We always take $a=0$ if $a\in M$.
Since $V_{\ol0},V_{\ol1}$ are $\Vir[M]$-modules such that
$\dim (V_\si)_\l\le1,\,\f\si\in\Z/2\Z,\l\in\F$, they have the form
$A_{a,b}(M),A_a(M),B_a(M)$ or their subquotients or the direct sum of the
two composition factors of $A_{0,0}(M)$.
Since we can interchange $V_{\ol0}$ and $V_{\ol1}$ if necessary,
it suffices to consider the following cases.
\par
{\it Case 1:} Suppose $V_{\ol0}=\F v_0$.
\par
Let $w_\nu=G_\nu v_0,\,\nu\in\a+M$, then applying $L_\mu,G_\eta$ to
$G_\nu v_0$, by (5.17), we have
$$
\matrix{
L_\mu w_\nu=L_\mu G_\nu v_0=
(\nu-{\mu\over2})w_{\mu+\nu},\,\f\mu\in M,\ \nu\in\a+M,\vs{2pt}\hfill\cr
G_\eta w_\nu=G_\eta G_\nu v_0=-G_\nu w_\eta,\,\f\eta,\nu\in \a+M.\hfill\cr
}
\eqno(5.20)$$
Since $G_\eta w_\nu=0$ for all $\eta+\nu\ne0$, we have
$$
(\nu-\eta)w_{2\eta+\nu}=L_{2\eta}w_\nu=G^2_\eta w_\nu=0.
\eqno(5.21)$$
Thus $w_\nu=0$ for all $\nu\in\a+M$ and $V$ is a trivial module.
\par
{\it Case 2:} Suppose $V_{\ol0}=A'_{0,0}(M)$.
\par
Then we can choose a basis $\{v_\mu\,|\,\mu\in M'\}$, where
$M'=M\setminus\{0\}$ such that
$$
L_\l v_{\mu}=(\l+\mu)v_{\l+\mu},\ \f\l\in M,\mu\in M'.
\eqno(5.22)$$
First fix $0\ne\mu_0\in M$. For any $\nu\in\a+M$, then $\nu-\mu_0\in\a+M$ and
we set $w'_\nu=G_{\nu-\mu_0}v_{\mu_0}$.
If $2\nu-\mu_0\ne0$, then
we have $G_{\nu-\mu_0}w'_\nu=L_{2(\nu-\mu_0)}v_{\mu_0}=
(2\nu-\mu_0)v_{2\nu-\mu_0}\ne0$ and so $w'_\nu\ne0$. Also
if $\nu_0\in\a+M$ such that $2\nu_0-\mu_0=0$, then
$G_{\nu_0}w'_{2\nu_0}=G^2_{\nu_0}v_{\mu_0}=L_{2\nu_0}v_{\mu_0}=
3\mu_0v_{2\mu_0}\ne0$. This shows that $w'_\nu\ne0$ for all $\nu\in\a+M$.
Therefore we can choose a basis $\{w_\nu\,|\,\nu\in\a+M\}$ of $V_{\ol1}$
such that there exists $b\in\F$ and
$$
L_\mu w_\nu=(\nu+b \mu)w_{\mu+\nu},
\eqno(5.23)$$
holds for $\mu\in M,\,\nu\in\a+M,\ \nu\ne0,\mu+\nu\ne0.$
For any $0\ne\nu\in\a+M$, applying $L_{2\nu}$ to $w_{-\nu}$, we have
$\nu(2b-1)w_\nu=L_{2\nu}w_{-\nu}=
G_\nu\cdot G_\nu w_{-\nu}\in G_\nu (V_{\ol0})_0=\{0\}.$
Thus $b={1\over2}$ and so (5.23) holds for all $\mu\in M,\,\nu\in\a+M$.
Now suppose
$$
G_\nu v_\mu=a_{\nu,\mu}w_{\nu+\mu},\
G_\nu w_\eta=b_{\nu,\eta}x_{\nu+\eta},\ 0\ne\mu\in M,\ \nu,\eta\in\a+M,
\eqno(5.24)$$
for some $a_{\nu,\mu},b_{\nu,\eta}\in\F$. We can suppose $a_{\mu,\nu}=1$
for some $\mu,\nu$. Apply $G_\nu,L_{2\nu}$ to (5.24), by (5.22-23),
we obtain
$$
\matrix{
a_{\nu,\mu}b_{\nu,\nu+\mu}=\mu+2\nu,\hfill&
b_{\nu,\eta}a_{\nu,\nu+\eta}=\eta+\nu,\vs{2pt}\hfill\cr
(2\nu+\mu)a_{\nu,\mu}=(2\nu+\mu)a_{\nu,\mu+2\nu},\hfill&
(3\nu+\eta)b_{\nu,\eta}=(\eta+\nu)b_{\nu,\eta+2\nu}.\hfill\cr
}
\eqno(5.25)$$
Now one can proceed exactly as in [S1, \S2] to prove
that $a_{\nu,\mu}=1$ and $b_{\nu,\eta}=\eta+\nu$, so that $V$ is the
quotient module of $\SA_0(M)$.
\par
{\it Case 3:} Suppose $V_{\ol0}=A_a(M),B_a(M),A_{0,0}(M)$ or
direct sum of two submodules $A'_{0,0}(M)\oplus\F v_0$.
\par
Similar to the proof of Cases 2, it is not difficult to show $V=\SA_a(M),\SB_a(M),
\SA_{0,0}(M)$ or direct sum of two submodules
$\SA'_{0,0}(M)\oplus\F v_0$
(this can not occur since $V$ is \nob{indecomposable}).
\par
{\it Case 4:} Suppose both $V_{\ol0}$ and $V_{\ol1}$ are simple
$\Vir[M]$-module of type $A_{a,b}(M)$.
\par
As in case 2, we can deduce similar equation as (5.25) and
the proof is exactly analogous to that of [S, Theorem 3.1]. We
obtain that $V$ is of type $\SA_{a,b}(M)$.
\qed
\par\
\par\
\par\small
\parskip 0.06 truein\baselineskip 3.5pt\lineskip 3.5pt
\cl{\bf References}
\par\ni\hi4ex\ha1
[DZ1]  D.Z. Djokovic and K. Zhao, Some simple subalgebras
of generalized Block algebras, {\it J. Alg.},
Vol.192, 74-101(1997).
\par\ni\hi4ex\ha1
[DZ2]   D.Z. Djokovic and K. Zhao, Derivations, isomorphisms, and
second cohomology of \nob{generalized} Witt algebras, {\it Trans. of
Amer. Math. Soc.}, Vol.350(2), 643-664(1998).
\par\ni\hi4ex\ha1
[DZ3]  D.Z.~Djokovic~and~K.~Zhao,~Some~infinite-dimensional~simple
Lie~algebras~in~\nob{characteristic} $0$ related to those of Block,
{\it J. Pure and Applied Algebras}, Vol.127(2), 153-165(1998).
\par\ni\hi4ex\ha1
[DZ4]  D.Z. Djokovic and K. Zhao, Generalized Cartan type $W$ Lie algebras
in characteristic zero, {\it J. Alg.},
Vol.195, 170-210(1997).
\par\ni\hi4ex\ha1
[DZ5]  D.Z. Djokovic and K. Zhao, Generalized Cartan type $S$ Lie algebras
in characteristic zero, {\it J. Alg.}, Vol.193, 144-179(1997).
\par\ni\hi4ex\ha1
[DZ6]  D.Z. Djokovic and K. Zhao, Second cohomology of Generalized
Cartan type $S$ Lie algebras in characteristic zero,
{\it J. Pure and Applied Algebras}, 136(1999), 101-126.
\par\ni\hi4ex\ha1
[KR] V. G. Kac and Raina, {\it Bombay lectures on highest weight
representations of infinite \nob{\it dimensional}  Lie algebras}, World Sci.,
Singapore, 1987.
\par\ni\hi4ex\ha1
[OZ1]  J.M. Osborn and K. Zhao, Infinite dimensional Lie algebras
of type L, {\it Comm. Alg.}, accepted for publication.
\par\ni\hi4ex\ha1
[OZ2]  J.M. Osborn and K. Zhao, Infinite dimensional Lie algebras
of generalized Block type, {\it Proc. of Amer. Math. Soc.},
Vol.127, No.6, 1641-1650(1999).
\par\ni\hi4ex\ha1
[OZ3]  J.M. Osborn and K. Zhao, Generalized Poisson bracket and Lie algebras
of type H in characteristic 0, {\it Math. Z.}, Vol.230, 107-143(1999).
\par\ni\hi4ex\ha1
[OZ4]  J.M. Osborn and K. Zhao, Generalized Cartan type K Lie algebras
in characteristic 0, {\it Comm. Alg.}, 25(1997), 3325-3360.
\par\ni\hi4ex\ha1
[OZ5]  J.M. Osborn and K. Zhao, Doubly $\bf Z$-graded Lie
algebras containing a Virasoro algebra, {\it J. Alg.}, 219(1999), 266-298.
\par\ni\hi4ex\ha1
[PZ] J. Patera and H. Zassenhaus, The higher rank Virasoro algebras,
{\it Comm. Math. Phys.} {\bf 136}(1991), 1-14.
\par\ni\hi4ex\ha1
[S] Y. C. Su, Harish-Chandra modules of the intermediate series over
the high rank Virasoro algebras and high rank super-Virasoro
algebras, {\it J. Math. Phys.} {\bf 35}(1994), 2013-2023.
\par\ni\hi4ex\ha1
[S1] Y. C. Su, Classification of Harish-Chandra modules over the
super-Virasoro algebras, {\it Comm. Alg.} {\bf 23}(1995), 3653-3675.
\end{document}